\documentclass[letterpaper, 10 pt, conference]{IEEEtran}
\usepackage{cite}
\usepackage{balance}

\usepackage{afterpage}
\usepackage{algorithm}
\usepackage[noend]{algpseudocode}
\usepackage{pifont}
\usepackage{wrapfig}
\usepackage{amsmath,amssymb}
\usepackage{bm}
\usepackage[authormarkup=none,final]{changes}
\usepackage{epigraph}
\usepackage{booktabs}
\usepackage{float}
\usepackage{graphicx}
\graphicspath{{Figures/}}
\usepackage{lettrine}
\usepackage{mathtools}
\usepackage{mathrsfs}
\usepackage{multirow}
\usepackage{hyperref}
\usepackage[percent]{overpic}
\usepackage[font={footnotesize,it}]{caption}
\usepackage{subcaption}
\usepackage{tikz}
\usepackage{upgreek}
\usepackage{xcolor}
\usepackage{xfrac}
\usepackage{siunitx}
\def\BState{\State\hskip-\ALG@thistlm}
\usepackage{enumitem}
\usepackage{blindtext}
\SetLabelAlign{parright}{\smash{\parbox[t]{\labelwidth}{\raggedleft#1}}}

\usepackage[left=48pt, right=48pt, bottom=76pt, top=50pt]{geometry}
\def\BibTeX{{\rm B\kern-.05em{\sc i\kern-.025em b}\kern-.08em
		T\kern-.1667em\lower.7ex\hbox{E}\kern-.125emX}}

\IEEEoverridecommandlockouts

\newtheorem{definition}{\textbf{Definition}}

\newtheorem{theorem}{\textbf{Theorem}}
\newtheorem{lemma}{\textbf{Lemma}}

\newenvironment{proof}{{\bf \emph{Proof.} }}{\hfill $\Box$ \\}

\newtheorem{proposition}{\textbf{Proposition}}
\newtheorem{corollary}{\textbf{Corollary}}
\newenvironment{proofof}[1]{{\bf \emph{Proof of #1.} }}{\hfill $\Box$ \par\smallskip}
\newcommand{\R}{\mathbb{R}}
\newcommand{\N}{\mathbb{N}}
\newcommand{\T}{\top}
\newcommand{\Lag}{\mathcal{L}}
\newcommand{\cV}{\check{\mathcal{V}}}
\newcommand{\cX}{\check{\mathcal{X}}}
\newcommand{\cQ}{\check{\mathcal{Q}}}
\newcommand{\Vc}{\check{V}}
\newcommand{\vc}{\check{v}}
\newcommand{\Wc}{\check{W}}
\newcommand{\wc}{\check{w}}
\newcommand{\Qc}{\check{Q}}
\newcommand{\qc}{\check{q}}
\newcommand{\Pc}{\check{P}}
\newcommand{\pc}{\check{p}}

\newcommand{\ext}{\operatorname*{ext}}
\newcommand{\bul}{\bullet}

\begin{document}

	\title{\LARGE Costate-Conditioned Dynamic Programming for Parallel-in-Time Solution of
	Equality-Constrained Optimal Control Problems}
	\author{ Tom Lefebvre${}^\star$
		\thanks{${}^\star$\{name.surname\}@ugent.be corresponding author
		}
		\thanks{
			The author is with the Dynamical Systems \& Control (DySC) research group, Department of Electromechanical, Systems and Metals Engineering, Faculty of Engineering and Architecture, Ghent University, Ghent, Belgium.
		}
	}

	\maketitle
	\thispagestyle{plain}
	\pagestyle{plain}

	\begin{abstract}
Gradient based trajectory optimisation algorithms spend most of their per-iteration time on a structured Newton--KKT system. The system corresponds to an equality-constrained linear-quadratic optimal control problem and dedicated Riccati-type solvers exhibit linear complexity in the horizon length. Parallel-in-time algorithms exhibit the potential to reduce the complexity to logarithmic in the horizon by distributing several subtrajectories over multiple cores and then stitch together the results by means of an associative scan. Existing exact formulations condition each subtrajectory on its terminal state. Each subproblem then carries the reachability to that terminal state explicitly as constraint rows. The associated reachability decision however is ill-conditioned for weakly actuated systems. In this work we condition each subtrajectory on the terminal costate instead. The resulting subproblems can now be computed by a free-end Riccati recursion without reachability carrying rows. This approach effectively compensates for the ill-conditioning of weakly actuated systems but requires additional bookkeeping with Lagrangian multipliers when stitching together the various subproblems. Numerically the method reproduces a direct solution of the  primal-dual Newton-KKT system to $10^{-12}$ in relative terms, on a weakly controllable chain of masses. The accurracy rivals that of a sequential recursion, whereas terminal-state conditioning loses up to six digits or fails. Operation counts on a problem with 64 subtrajectories and a horizon of 4096 stages gives a speed-up of 26 compared with a sequential solve.
	\end{abstract}

\section{Introduction}\label{sec:intro}

Direct multiple shooting \cite{bock1984multiple} transcribes a trajectory optimisation problem into a
nonlinear program with stagewise structure. Interior-point methods of the IPOPT type
\cite{waechter2006implementation} solve such programs robustly. Their per-iteration cost is dominated by
the solution of a Newton--KKT system. After elimination of slack variables and inequality multipliers, which are
local to each stage, the solution to this system is equivalent to the solution of an equality-constrained
linear-quadratic (EC-LQ) optimal control problem \cite{rao1998application}. For unconstrained problems this solution can be obtained efficiently using a Ricatti recursion. FATROP \cite{vanroye2023fatrop}
tackles equality constrained problem using a generalised Riccati recursion that handles stagewise equality constraints exactly
\cite{vanroye2024generalization,sideris2010riccati,laine2019efficient}.
By construction, a sequential Ricatti recursion has time complexity $\mathcal{O}(T)$ in the horizon length, $T$.

Parallel-in-time methods have the potential to exploit modern parallel hardware such as GPUs. Theoretically they allow logarithmic complexity, that is $\mathcal{O}(\log(T))$, which may result in significant practical speed-up, especially for long horizon problems. Consequently, several works have produced parallel-in-time methods tailored to dynamic programming problems. Wright \cite{wright1991partitioned}
partitioned the horizon for parallel dynamic programming. Nielsen and Axehill
\cite{nielsen2014ologn,nielsen2019direct} computed Newton steps for model predictive control in
$\mathcal{O}(\log T)$ parallel steps by reducing parametric subproblems to a master problem of the same form.
S\"arkk\"a and Garc\'{\i}a-Fern\'andez \cite{sarkka2023temporal} formulated the unconstrained LQ problem as an
associative scan \cite{blelloch1989scans} over conditional value functions in the sense of Calvet and Viargues
\cite{calvet1995invariant,houska2017costtotravel}. Recent GPU based trajectory optimization solvers mostly adopted the associative scan approach
\cite{iacob2025parallel,amatucci2026primaldual}.

Stagewise equality constraints, arising from e.g. contact conditions, waypoints or physical invariants, have been added to
parallel-in-time methods in two ways. The first regularises the multipliers of the equality constraints,
which allows them to be eliminated stage by stage before parallel solution \cite{jallet2024parallel,sousapinto2025dual}. In such methods every Newton step satisfies the linearised constraints
only up to the regularisation and feasibility is reached through the multiplier updates. The second approach ensures exact constraint satisfaction \cite{laine2019parallelizing,lefebvre2026temporal}. Here the constrained problem is treated similarly to the approach of S\"arkk\"a, introducing terminal state conditioned value functions for every subtrajectory. To account for the equality constraints, the conditional value function is appended with an implicit domain function tracking the reachability of the conditioning state. In a linear-quadratic setting, recombination of the independent subproblems may become ill-conditioned when the subtrajectories are short or when the system is weakly actuated. 

In this work we build on the exact approach but condition each subtrajectory on the terminal costate. By construction this avoids ill-conditioning issues when recombining the independent subproblems. The disadvantage is that recombination now requires additional bookkeeping. To that end we propose a two step recombination procedure that preserves the associativity and thus supports parallel-in-time execution. We further discuss how to recover the complete primal-dual solution and provide some numerical case studies to evaluate accuracy and practical complexity.

To keep the paper self-contained, Section~\ref{sec:statecond} restates the state-conditioned formulation of
\cite{lefebvre2026temporal} in its notation, which is used throughout. All proofs and the explicit coefficient
formulas of the LQ specialisation are given in the appendices.

\section{Problem Statement}\label{sec:problem}

\subsection{Notation} We write $\alpha_{t_1:t_2}=\{\alpha_{t_1},\dots,\alpha_{t_2}\}$ for the collection of a variable or function $\alpha$
over the discrete times $t_1,\dots,t_2$, and $\N_{a:b}=\{a,\dots,b\}$. All parametric problems considered in this work are equality-constrained quadratic programs. With the $\min$ operator we denote their stationary value. With $\ext$ we denote stationary with respect to a multiplier-typ variable. Note that any stationarity condition that does not involve the variable becomes a constraint row of the result.

\subsection{Equality-constrained optimal control}

Consider the equality-constrained,
deterministic, discrete-time optimal control problem
\begin{equation}\label{eq:ocp}
\begin{aligned}
\min_{u_{0:T-1}}\ & \sum\nolimits_{t=0}^{T-1}c_t(x_t,u_t)+c_T(x_T)\\
\text{s.t. }\ & 0=x_{t+1}-f_t(x_t,u_t), t\in\N_{0:T-1},\\
& 0=g_t(x_t,u_t), t\in\N_{0:T-1},\\
& 0=g_T(x_T),
\end{aligned}
\end{equation}
with state, $x_t\in\R^{n_x}$, input, $u_t\in\R^{n_u}$, stage cost, $c_t$, terminal cost, $c_T$, dynamics, $f_t$, stage
constraint, $g_t$, with $n_{g,t}$ rows, terminal constraint, $g_T$, with $n_{g,T}$ rows, and $x_0$ given.

In the linear-quadratic (LQ) setting these functions adopt the following form
\begin{equation}\label{eq:lqdata}
\begin{aligned}
c_t(x_t,u_t)&=\tfrac12\begin{pmatrix}1\\x_t\\u_t\end{pmatrix}^{\!\T}
\begin{pmatrix}2c_{t,0}&\bul&\bul\\c_{t,x}&C_{t,xx}&\bul\\c_{t,u}&C_{t,ux}&C_{t,uu}\end{pmatrix}
\begin{pmatrix}1\\x_t\\u_t\end{pmatrix},\\
c_T(x_T)&=\tfrac12\begin{pmatrix}1\\x_T\end{pmatrix}^{\!\T}
\begin{pmatrix}2c_{T,0}&\bul\\c_{T,x}&C_{T,xx}\end{pmatrix}
\begin{pmatrix}1\\x_T\end{pmatrix},\\
f_t(x_t,u_t)&=F_{t,x}x_t+F_{t,u}u_t+f_{t,0},\\
g_t(x_t,u_t)&=G_{t,x}x_t+G_{t,u}u_t+g_{t,0},\\
g_T(x_T)&=G_{T,x}x_T+g_{T,0},
\end{aligned}
\end{equation}
where $\bul$ denotes entries implied by symmetry. 

\subsection{Primal-dual solution}\label{sec:ipm}

The Lagrangian corresponding to \eqref{eq:ocp} is defined as
\begin{equation}\label{eq:lagrangian}
\begin{multlined}
\Lag={}\sum\nolimits_{t=0}^{T-1}\Big(c_t(x_t,u_t)+\gamma_t^\T g_t(x_t,u_t) + \\\pi_{t+1}^\T\big(f_t(x_t,u_t)-x_{t+1}\big)\Big) +c_T(x_T)+\gamma_T^\T g_T(x_T),
\end{multlined}
\end{equation}
with multipliers $\pi_{1:T}$\footnote{The given sign convention ensures that $\pi_t$ refers to the costate, that is the gradient at $x_t$ of the optimal cost-to-go, including its constraint terms.} and $\gamma_{0:T}$. We assume that the problem is well posed and a unique optimal solution exists.

An IPOPT-type primal--dual interior-point method \cite{waechter2006implementation} applied to a
multiple-shooting nonlinear program computes a search direction from the primal--dual system. After the
stagewise elimination of slacks and inequality multipliers, this system corresponds with the optimality system of an LQ
problem \cite{rao1998application,vanroye2023fatrop}. The
linear solver must therefore return the complete solution including the primal step, $x_{0:T}$ and $u_{0:T-1}$, and, the dual step consisting of all multipliers $\pi_{1:T}$ and $\gamma_{0:T}$. 

\section{State-conditioned Parallel Programming}\label{sec:statecond}

First we detail the solution proposed in \cite{laine2019efficient,lefebvre2026temporal}.

\subsection{Value tuples}

The \emph{value function} $V_t(x_t)$ is the optimal cost of \eqref{eq:ocp} restricted to the stages $t,\dots,T$ and
started at $x_t$. The value function is defined on the \emph{value domain} $\mathcal{X}_t$, which is the set of states from which all future
constraints can be satisfied, represented implicitly by the \emph{value domain function} $w_t$. The domain is defined as $\mathcal{X}_t=\{x_t\in\R^{n_x}:0=w_t(x_t)\}$.

A \emph{tuple} consists of a cost element and a constraint element that share their argument. Correspondingly we define the \emph{value
tuple}, $\mathcal{V}_t$, and the \emph{stage tuple}, $\mathcal{C}_t$, as
\begin{equation}
\begin{aligned}
\mathcal{V}_t(x_t)&=(V_t(x_t),w_t(x_t)),\label{eq:vtuple}\\
\mathcal{C}_t(x_t,u_t)&=(c_t(x_t,u_t),g_t(x_t,u_t)).
\end{aligned}
\end{equation}

For tuples $\mathcal{A}=(a,f)$ and $\mathcal{B}=(b,g)$ the operator $\oplus$ adds the cost elements and stacks the
constraint elements, $\mathcal{A}\oplus\mathcal{B}=(a+b,\{f,g\})$. Minimising a tuple means minimising its cost element
subject to its constraint elements. The result is again a tuple whose constraints collect the conditions on
the remaining arguments for which the minimisation is feasible.

Similar to regular value functions, these value tuples allow a backward recursive evaluation.

\begin{lemma}[{\cite[Lemma~1]{lefebvre2026temporal}}]\label{lem:bellman}
With the $Q$-tuple $\mathcal{Q}_t(x_t,u_t)=\mathcal{C}_t(x_t,u_t)\oplus\mathcal{V}_{t+1}(f_t(x_t,u_t))$,
\begin{equation}
\begin{aligned}
\mathcal{V}_t(x_t)&=\min_{u_t}\ \mathcal{Q}_t(x_t,u_t),\label{eq:bellman}\\
u^*_t(x_t)&=\arg\min_{u_t}\ \mathcal{Q}_t(x_t,u_t),
\end{aligned}
\end{equation}
with boundary condition $\mathcal{V}_T(x_T)=(c_T(x_T),g_T(x_T))$.
\end{lemma}

\subsection{State-conditioned value tuples}

The state-conditioned formulation parametrises a subtrajectory $[t,t']$ by its initial and its terminal state.
The \emph{state-conditioned value function}
\begin{equation}\label{eq:cond}
\begin{aligned}
&V_{t\to t'}(x_t,x_{t'})=\min_{u_{t:t'-1}}\ \sum\nolimits_{\tau=t}^{t'-1}c_\tau(x_\tau,u_\tau)\\
&\quad\text{s.t. }0=x_{\tau+1}-f_\tau(x_\tau,u_\tau),\ \tau\in\N_{t:t'-1},\\
&\quad\phantom{\text{s.t. }}0=g_\tau(x_\tau,u_\tau),\ \tau\in\N_{t:t'-1},\\
&\quad\phantom{\text{s.t. }}0=x_{t'}-x_\tau,\ \tau=t',
\end{aligned}
\end{equation}
is defined on the state-conditioned value domain $\mathcal{X}_{t\to t'}=\{(x_t,x_{t'}):0=w_{t\to t'}(x_t,x_{t'})\}$. The
\emph{state-conditioned value tuple} and $Q$-tuple are defined as
\begin{equation}
\begin{aligned}
\mathcal{V}_{t\to t'}(x_t,x_{t'})&=(V_{t\to t'}(x_t,x_{t'}),w_{t\to t'}(x_t,x_{t'})),\\
\mathcal{Q}_{t\to t'}(x_t,u_t,x_{t'})&=\mathcal{C}_t(x_t,u_t)\oplus\mathcal{V}_{t+1\to t'}(f_t(x_t,u_t),x_{t'}).
\end{aligned}
\end{equation}

State-conditioned tuples again allow a backward recursive evaluation. 
\begin{theorem}[{\cite[Theorem~1]{lefebvre2026temporal}}]\label{thm:condrec}
The state-conditioned value tuple and optimal policy satisfy
\begin{equation}
\begin{aligned}
\mathcal{V}_{t\to t'}(x_t,x_{t'})&=\min_{u_t}\ \mathcal{Q}_{t\to t'}(x_t,u_t,x_{t'}),\\
u^*_{t\to t'}(x_t,x_{t'})&=\arg\min_{u_t}\ \mathcal{Q}_{t\to t'}(x_t,u_t,x_{t'}),
\end{aligned}
\end{equation}
with the boundary condition
\begin{equation}
\begin{aligned}
V_{t'\to t'}(x,x')&=0,\label{eq:lvinit1}\\
w_{t'\to t'}(x,x')&=x-x'
\end{aligned}
\end{equation}
\end{theorem}

Interestingly, state-conditioned value tuples can be combined to form another state-conditioned value tuple that now spans the mutual time interval. 



\begin{theorem}[{\cite[Theorem~2]{lefebvre2026temporal}}]\label{thm:lvcomb}
For $t<\tau<t'\le T$,
\begin{align}
\mathcal{V}_{t\to t'}(x_t,x_{t'})&=\min_{x_\tau}\ \mathcal{V}_{t\to\tau}(x_t,x_\tau)\oplus\mathcal{V}_{\tau\to t'}(x_\tau,x_{t'}), \notag\\
x^*_{\tau,t\to t'}(x_t,x_{t'})&=\arg\min_{x_\tau}\ \mathcal{V}_{t\to\tau}(x_t,x_\tau)\oplus\mathcal{V}_{\tau\to t'}(x_\tau,x_{t'}), \notag\\
\mathcal{V}_t(x_t)&=\min_{x_\tau}\ \mathcal{V}_{t\to\tau}(x_t,x_\tau)\oplus\mathcal{V}_\tau(x_\tau), \notag\\
x^*_{\tau,t}(x_t)&=\arg\min_{x_\tau}\ \mathcal{V}_{t\to\tau}(x_t,x_\tau)\oplus\mathcal{V}_\tau(x_\tau), \label{eq:lvc2}
\end{align}
\end{theorem}

The combination rule above is \textit{associative}
\cite[Corollary~4]{lefebvre2026temporal} and supports parallel-in-time solution of the problem.

\subsection{Parallel solution}\label{sec:parsolstatecond}

Assume the horizon of \eqref{eq:ocp} is cut at time instants $t_0=0<t_1<\dots<t_{N-1}<t_N=T$. 
\begin{itemize}
    \item In an \emph{initial pass} the tuples $\mathcal{V}_{t_i\to t_{i+1}}$, $i\in\N_{0:N-2}$, and $\mathcal{V}_{t_{N-1}}$ are computed independently on different cores.
    \item  In an \emph{upward sweep}, the tuples are combined, two by two, applying Theorem~\ref{thm:lvcomb}, therewith spanning a binary tree in $\lceil\log_2N\rceil$ levels. Every node in the tree only needs to store its \textit{connecting-state map}, $x^*_{\tau,t\to t'}(x_t,x_{t'})$. 
    \item In the
\emph{downward sweep} every node evaluates its map with the arguments received from its parent and passes,
$(x_t^*,x^*_\tau)$ upstream, and, $(x^*_\tau,x_{t'}^*)$ downstream.
\item  A \emph{final pass} simulates every
subtrajectory in closed loop with $u^*_{t\to t_{i+1}}(x_t,x^*_{t_{i+1}})$, $t\in[t_i,t_{i+1}]$ producing the full (primal) solution, consisting of $x_{0:T}^*$ and $u_{0:T-1}^*$.
\end{itemize}

 The complexity of the total procedure is $\mathcal{O}(T/N+\log_2N)$.

\section{Parametric quadratic programs}\label{sec:pqp}

In the LQ setting every decision from section~\ref{sec:statecond} reduces to a parametric quadratic program (p-QP) of the form
\begin{equation}\label{eq:pqp}
\min_u\ a(x,u) \text{ s.t. }0=b(x,u)=B_xx+B_uu+b_0,
\end{equation}
with \emph{parameter}, $x$, optimization \emph{variable}, $u$, a quadratic cost, $a$ (defined as in \eqref{eq:lqdata}) and an affine constraint, $b$. 

This section restates the null-space method of \cite{lefebvre2026temporal} and elaborates on the associated multipliers which are essential to understand the additional bookkeeping that is required at a cut on which the rest of the paper relies.

\subsection{Solution}

Consider an orthogonal transformation $T=[\tilde T\ \hat T]$ of the constraint, separating the rows that $u$ can
satisfy from those it cannot
\begin{equation}\label{eq:pqpT}
T^\T B_u=\begin{pmatrix}\tilde T^\T B_u\\\hat T^\T B_u\end{pmatrix}=\begin{pmatrix}\tilde B_u\\0\end{pmatrix},
\end{equation}
with $\tilde B_u$ of full row rank. The rows along $\tilde T$ are \emph{used}, i.e. $u$ satisfies them, and its remaining
freedom is spent on the cost. The rows along $\hat T$ are \emph{transferred upstream}, i.e. they do not involve $u$, so they
become a condition on the parameter, $x$, i.e. $0=\hat B_xx+\hat b_0$ with $\hat B_x=\hat T^\T B_x$ and $\hat b_0=\hat T^\T b_0$.

The solution of \eqref{eq:pqp} is an affine policy and a new tuple in the parameter, $x$. The rows from the constraint element of the
solution tuple, $\hat{b}(x)$, are handed over to the upstream problem that decides the parameter, $x$.
\begin{equation}
\begin{aligned}
u^*(x)&=U_xx+u_0,\label{eq:pqppol}\\
\big(\hat a(x),\hat b(x)\big)&=\big(a(x,u^*(x)),\hat B_xx+\hat b_0\big),
\end{aligned}
\end{equation}

\subsection{Lagrangian multipliers}\label{sec:pqpmult}

The Lagrangian of the p-QP \eqref{eq:pqp} is given by $\mathcal{L} = a(x,u)+\eta^\T b(x,u)$. Let us first decompose the multiplier with the same
transformation used to transform the constraint,
\begin{equation}\label{eq:pqpmult1}
\eta=\tilde T\tilde\eta+\hat T\hat\eta.
\end{equation}

\begin{proposition}\label{prop:pqp}
Let $x$ satisfy $\hat b(x)=0$. 
\begin{itemize}
    \item Stationarity in $u$ reads
\begin{equation}\label{eq:pqpmult2}
\nabla_ua(x,u^*(x))+\tilde B_u^\T\tilde\eta=0,
\end{equation}
which determines $\tilde\eta$ uniquely and does not involve $\hat\eta$.
    \item For every $\hat\eta$, the derivative of the
Lagrangian, $\mathcal{L},$ with respect to the parameter at the solution is
\begin{equation}\label{eq:pqpsens}
\nabla_xa(x,u^*(x))+B_x^\T\eta=\nabla_x\hat a(x)+\hat B_x^\T\hat\eta.
\end{equation}
\end{itemize}
\end{proposition}

Proposition~\ref{prop:pqp} is the principle that will come to organise and recover all multipliers in this paper.
\begin{itemize}[leftmargin=*]
\item The multiplier of the \emph{used} rows, $\tilde\eta$, is determined by the p-QP itself, as soon as its
parameter is known.
\item The multiplier of the \emph{transferred} rows, $\hat\eta$, is \emph{not} determined by the p-QP. The transferred
rows constrain the parameter, so their multiplier is determined by the decision that determines the parameters whilst satisfying these rows.
\item The \emph{price of the parameter} is \eqref{eq:pqpsens} and consists of the gradient of the optimal value, in which
the used rows are priced implicitly, and of the transferred rows priced at $\hat\eta$. The policy, $u^*(x)$, does not depend on
$\hat\eta$, but the price does. A decision that ignores $\hat\eta$ obtains the correct $u^*$ but the wrong price and thus ultimately the wrong multipliers.
\end{itemize}



\subsection{Application to LQ state-conditioned tuples}\label{sec:lvlq}

Applied to the setting \eqref{eq:ocp}-\eqref{eq:lqdata}, the state-conditioned value tuple is
\begin{equation}
\begin{aligned}
V_{\star}={}&
\tfrac12\begin{pmatrix}1\\x_t\\x_{t'}\end{pmatrix}^{\!\T}
\begin{pmatrix}2v_{\star,0}&\bul&\bul\\v_{\star,x}&V_{\star,xx}&\bul\\v_{\star,x'}&V_{\star,x'x}&V_{\star,x'x'}\end{pmatrix}
\begin{pmatrix}1\\x_t\\x_{t'}\end{pmatrix} 
\\
w_{\star}={}&W_{t\to t',x}x_t+W_{t\to t',x'}x_{t'}+w_{t\to t',0},\label{eq:lvw}
\end{aligned}
\end{equation}
where arguments are omitted for brevity and $\star=t\to t'$.

The recursion from Theorem~\ref{thm:condrec} applied to \eqref{eq:ocp}-\eqref{eq:lqdata} reduces to a p-QP with variable $u_t$ and parameter $(x_t,x_{t'})$. The (terminal) state-conditioned optimal policy is
\begin{equation}\label{eq:lvpol}
u^*_{t\to t'}(x_t,x_{t'})=U_{t\to t',x}x_t+U_{t\to t',x'}x_{t'}+u_{t\to t',0}.
\end{equation}

Similarly, the recombination from Theorem \ref{thm:lvcomb}
reduces to a p-QP with variable $x_\tau$ and parameter $(x_t,x_{t'})$, and its policy is the connecting-state map
\begin{equation}\label{eq:lvX}
x^*_{\tau,t\to t'}(x_t,x_{t'})=X_{\tau,t\to t',x}x_t+X_{\tau,t\to t',x'}x_{t'}+x_{\tau,t\to t',0}.
\end{equation}

This approach comes with two limitations. The boundary condition in \eqref{eq:lvinit1} has $n_x$ rows. Each step of the
recursion uses at most $n_u$ of them and transfers the others upstream. So for at least $\lceil n_x/n_u\rceil$ stages every p-QP
carries up to $n_x+n_{g,t}$ rows and decides which terminal states are reachable. Where they are, $V_{t\to t',x'x'}$ behaves
like the inverse of the reachability Gramian of the subtrajectory, which is nearly singular for short subtrajectories
and weakly actuated systems. The associated decisions then become ill-conditioned. 

\section{Costate-conditioned Parallel Programming}\label{sec:costate}

To remedy the ill-conditioning and multiple recover problems emerging with state-conditioned value tuples, we propose to condition the tuples on the \textit{terminal costate}. 

\subsection{Costate-conditioned value tuples}

\begin{definition}\label{def:costate}
For $t<t'\le T$ and $\lambda_{t'}\in\R^{n_x}$, the \emph{costate-conditioned value function} is
\begin{equation}\label{eq:defcostate}
\begin{aligned}
&\Vc_{t\to t'}(x_t,\lambda_{t'})=\min_{u_{t:t'-1}}\ \sum\nolimits_{\tau=t}^{t'-1}c_\tau(x_\tau,u_\tau)+\lambda_{t'}^\T x_{t'}\\
&\quad\text{s.t. }0=x_{\tau+1}-f_\tau(x_\tau,u_\tau), \tau\in\N_{t:t'-1},\\
&\quad\phantom{\text{s.t. }}0=g_\tau(x_\tau,u_\tau), \tau\in\N_{t:t'-1},
\end{aligned}
\end{equation}
defined on the domain $\cX_{t\to t'}=\{x_t:0=\wc_{t\to t'}(x_t)\}$ of initial states from which the stage
constraints can be satisfied. 

The \emph{costate-conditioned value tuple} is 
\begin{equation}
\begin{aligned}
\cV_{t\to t'}(x_t,\lambda_{t'})&=(\Vc_{t\to t'}(x_t,\lambda_{t'}),\wc_{t\to t'}(x_t)).
\end{aligned}
\end{equation}
\end{definition}

Compared with \eqref{eq:cond}, the \textit{explicit} terminal condition is removed and a linear terminal cost, $\lambda^\top_{t'} x_{t'}$, is added.  Since by the envelope theorem we have that $x_{t'} = \nabla_{\lambda_{t'}^*} \Vc_{t\to t'}(x_t,\lambda_t)$, the multiplier thus acts as an \textit{implicit} terminal state configuring parameter. We will come back to this insight later.


Once again, these tuples allow a recursive evaluation.

\begin{theorem}\label{thm:leaf}
With the $Q$-tuple $\cQ_{t\to t'}(x_t,u_t,\lambda_{t'})=\mathcal{C}_t(x_t,u_t)\oplus\cV_{t+1\to t'}(f_t(x_t,u_t),\lambda_{t'})$, the costate-conditioned value tuple and optimal policy satisfy
\begin{equation}
\begin{aligned}
\cV_{t\to t'}(x_t,\lambda_{t'})&=\min_{u_t}\ \cQ_{t\to t'}(x_t,u_t,\lambda_{t'}),\label{eq:leafrec}\\
u^*_{t\to t'}(x_t,\lambda_{t'})&=\arg\min_{u_t}\ \cQ_{t\to t'}(x_t,u_t,\lambda_{t'}),
\end{aligned}
\end{equation}
with boundary condition $\cV_{t'\to t'}(x,\lambda)=(\lambda^\T x,\emptyset)$.
\end{theorem}

\subsection{LQ costate-conditioned value tuples}\label{sec:costatelq}

Applied again to the setting \eqref{eq:ocp}-\eqref{eq:lqdata}, the costate-conditioned value tuple and policy are
\begin{equation}
\begin{aligned}
\Vc_{\star}={}&
\tfrac12\begin{pmatrix}1\\x_t\\\lambda_{t'}\end{pmatrix}^{\!\T}
\begin{pmatrix}2\vc_{\star,0}&\bul&\bul\\\vc_{\star,x}&\Vc_{\star,xx}&\bul\\\vc_{\star,\lambda}&\Vc_{\star,\lambda'x}&\Vc_{\star,\lambda\lambda}\end{pmatrix}
\begin{pmatrix}1\\x_t\\\lambda_{t'}\end{pmatrix} 
\label{eq:lqcostate}\\
\wc_{\star}={}&\Wc_{t\to t',x}x_t+\wc_{t\to t',0},
\end{aligned}
\end{equation}
Arguments are again omitted for brevity and $\star=t\to t'$. 

The corresponding costate-conditioned optimal policy is
\begin{equation}
    u^*_{t\to t'}={}U_{t\to t',x}x_t+U_{t\to t',\lambda}\lambda_{t'}+u_{t\to t',0},\label{eq:lqpolicy}
\end{equation}

Again, by the envelope theorem the terminal state is
\begin{equation}\label{eq:termstate}
\nabla_{\lambda_{t'}}\Vc_{t\to t'}=x^*_{t'}=\Vc_{t\to t',\lambda x}x_t+\vc_{t\to t',\lambda}+\Vc_{t\to t',\lambda\lambda}\lambda_{t'}.
\end{equation}

The recursion of Theorem \ref{thm:leaf} also determines a p-QP with variable, $u_t$, parameter, $(x_t,\lambda_{t'})$, and  boundary condition, $\Vc_{t'\to t',\lambda x}=I$. 

All coefficients defined here are given in the appendices. 

\begin{proposition}\label{prop:leaf} Consider the costate-conditioned tuple
\begin{itemize}
    \item The constraint element never depends on $\lambda_{t'}$.
    \item Parameters $\Vc_{t\to t',xx}$, $\vc_{t\to t',x}$, $\Wc_{t\to t',x}$, $\wc_{t\to t',0}$,
$U_{t\to t',x}$, $u_{t\to t',0}$ and all factorisations coincide with those of a conventional recursion on $[t,t']$ with zero
terminal cost and no terminal constraint
\end{itemize}
\end{proposition}

A tuple therefore never carries reachability rows. The only rank decisions concern the stage constraints, as in a sequential recursion. As a result, the relation between the initial and the terminal state, which the
state-conditioned tuple specifies with its domain, is now contained \textit{entirely} in the cost element. 

By \eqref{eq:termstate}, the
terminal states that the subtrajectory can reach from $x_t$ form the affine set
\begin{equation}\label{eq:reachset}
\mathcal{R}_{t\to t'}(x_t)=\Vc_{t\to t',\lambda x}x_t+\vc_{t\to t',\lambda}+\operatorname{range}\big(\Vc_{t\to t',\lambda\lambda}\big).
\end{equation}

Thus $\Vc_{t\to t',\lambda\lambda}$ adopts the role of a reachability Gramian. 

\subsection{Combination of LQ costate-conditioned tuples}\label{sec:comb}

Combining two costate-conditioned tuples requires additional bookkeeping. When two consecutive subtrajectories $[t,\tau]$ and $[\tau,t']$ are combined at $\tau$, the terminal
state, $x_\tau$, of the upstream subtrajectory must equal the initial state of the downstream substrajectory, also $x_\tau$, and, must be an element of the \emph{domain} of the downstream tuple. 

The downstream subtrajectory accepts initial states, $x_\tau$, in
\begin{equation}\label{eq:accset}
\mathcal{A}_{\tau\to t'}=\{x_\tau:0=\Wc_{\tau\to t',x}x_\tau+\wc_{\tau\to t',0}\}.
\end{equation}

The combination must therefore determine the overlap of
$\mathcal{R}_{t\to\tau}(x_t)$ and $\mathcal{A}_{\tau\to t'}$, and the corresponding connecting state, $x_\tau$. 
We cannot satisfy both requirements with the same multiplier. To that end, we dualise the \textit{downstream domain} with a new multiplier $\nu_{\tau}$ and \textit{continuity} with the multiplier
$\lambda_\tau$. The dual variable, $\nu_\tau$, is used to push the decision onto the overlap $\mathcal{R}_{t\to\tau}(x_t)\cap\mathcal{A}_{\tau\to t'}$. The effective terminal cost of the upstream subtrajectory then becomes
$\lambda_\tau^\T x_\tau+\nu_{\tau}^\T(\Wc_{\tau\to t',x}x_\tau+\wc_{\tau\to t',0})$. This cost is linear in $x_\tau$ with slope
\begin{equation}\label{eq:mueff}
\mu_\tau=\lambda_\tau+\Wc_{\tau\to t',x}^\T\nu_{\tau}.
\end{equation}

To absorb the downstream domain in the upstream problem, we can now simply substitute $\lambda_\tau \leftarrow \mu_\tau$ in \eqref{eq:defcostate}, add the linear cost, $\wc_{\tau\to t',0}^\top\nu_\tau$, substitute \eqref{eq:termstate} into $\wc_{\tau\to t'}$ and then solve for $\nu_\tau$. The resulting problem then corresponds with \eqref{eq:defcostate} augmented with the downstream domain, or equivalently, with \eqref{eq:defcostate} evaluated at the \emph{effective costate},
$\mu_\tau$, with value $\Vc_{t\to\tau}(x_t,\mu_\tau)+\nu_{\tau}^\T\wc_{\tau\to t',0}$. 

 The full combination determines $\nu_{\tau}$ first (absorption step), for given $\lambda_\tau$,
 and then $\lambda_\tau$ and $x_\tau$ (junction step).

\subsubsection{Absorption step} The absorption solves the problem
\begin{equation}
\label{eq:step1}
    \cV'_{t\to\tau}(x_t,\lambda_\tau) =\ext_{\nu_\tau} \cV_{t\to \tau}(x_t,\mu_\tau) \oplus \left(\wc_{\tau\to t',0}^\top\nu_\tau,\emptyset\right)
\end{equation}
The parameters are
$(x_t,\lambda_\tau)$ and the variable is $\nu_\tau$.

Substituting \eqref{eq:mueff} in \eqref{eq:termstate} and then \eqref{eq:termstate} into the domain function learns that a shift $\Wc_{\tau\to t',x}^\T\nu_{\tau}$ of the terminal, $\lambda_{t'}$ moves the downstream domain value by $\Wc_{\tau\to t',x}\Vc_{t\to\tau,\lambda\lambda}\Wc_{\tau\to t',x}^\T\nu_{\tau}$, so the matrix
$\Wc_{\tau\to t',x}\Vc_{t\to\tau,\lambda\lambda}$ plays the role of $B_u$ in \eqref{eq:pqpT} and \eqref{eq:step1} behaves as the dual of the p-QP \eqref{eq:pqp}.

Let $T_\tau=[\tilde T_\tau\ \hat T_\tau]$ be
so that $\hat T_\tau^\T\Wc_{\tau\to t',x}\Vc_{t\to\tau,\lambda\lambda}=0$ and $\tilde T_\tau^\T\Wc_{\tau\to t',x}\Vc_{t\to\tau,\lambda\lambda}$ of full
row rank, let $\tilde W_{\tau\to t',x}=\tilde T_\tau^\T\Wc_{\tau\to t',x}$ and $\hat W_{\tau\to t',x}=\hat T_\tau^\T\Wc_{\tau\to t',x}$, and decompose the
multiplier as in \eqref{eq:pqpmult1},
\begin{equation}
\nu_{\tau}=\tilde T_\tau\tilde\nu_\tau+\hat T_\tau\hat\nu_\tau.\label{eq:nusplit}
\end{equation}
\begin{itemize}[leftmargin=*]
\item The rows along $\tilde T_\tau$ are used. Stationarity in $\tilde\nu_\tau$ requires that its terminal state at the effective costate satisfies them, which
determines the optimal parametric multiplier
\begin{equation}\label{eq:numap}
\tilde\nu^*_{\tau,t\to t'}(x_t,\lambda_\tau)=N_{\tau,t\to\tau,x}x_t+N_{\tau,t\to\tau,\lambda}\lambda_\tau+n_{\tau,t\to\tau,0}.
\end{equation}
This is the price that moves the terminal state within $\mathcal{R}_{t\to\tau}(x_t)$ onto $\mathcal{R}_{t\to\tau}(x_t)\cap\mathcal{A}_{\tau\to t'}$, at
least cost.
\item The rows along $\hat T_\tau$ are transferred upstream. We have that $\hat W_{\tau\to t',x}x_\tau$ takes the same value
on all of $\mathcal{R}_{t\to\tau}(x_t)$. These rows result in transferred domain rows on $x_t$,
\begin{equation}\label{eq:newrows}
0=\hat W_{\tau\to t',x}\big(\Vc_{t\to\tau,\lambda x}x_t+\vc_{t\to\tau,\lambda}\big)+\hat T_\tau^\T\wc_{\tau\to t',0},
\end{equation}
which state that the overlap is non-empty. By Proposition~\ref{prop:pqp}, the associated multiplier $\hat\nu_\tau$ is
\emph{not} determined here but by the decision that fixes $x_t$ and thus by an ancestor in the tree.
\end{itemize}

Substituting \eqref{eq:numap} in \eqref{eq:step1} gives the \emph{absorbed tuple}, $\cV'_{t\to\tau}(x_t,\lambda_\tau)$, which corresponds with the costate-conditioned tuple
of $[t,\tau]$ with the additional terminal constraint $x_\tau\in\mathcal{A}_{\tau\to t'}$. The domain $\wc'_{t\to\tau}$ is
$\wc_{t\to\tau}$ stacked with \eqref{eq:newrows}.

\subsubsection{Junction step} After absorption of the downstream domain, only continuity couples the subtrajectories. The junction solves the p-QP
\begin{equation}\label{eq:c1}
\begin{multlined}
\cV_{t\to t'}(x_t,\lambda_{t'})=\min_{x_\tau}\ext_{\lambda_\tau}\ \cV'_{t\to\tau}(x_t,\lambda_\tau)
\dots \\\oplus\big(-\lambda_\tau^\T x_\tau+\Vc_{\tau\to t'}(x_\tau,\lambda_{t'}),\emptyset\big).
\end{multlined}
\end{equation}

The stationarity conditions in $\lambda_\tau$ and $x_\tau$ are
\begin{equation}
\begin{aligned}
x_\tau&=\Vc'_{t\to\tau,\lambda x}x_t+\vc'_{t\to\tau,\lambda}+\Vc'_{t\to\tau,\lambda\lambda}\lambda_\tau,\label{eq:stat1}\\
\lambda_\tau&=\Vc_{\tau\to t',xx}x_\tau+\vc_{\tau\to t',x}+\Vc_{\tau\to t',\lambda x}^\T\lambda_{t'}.
\end{aligned}
\end{equation}

The first demands continuity. The terminal state of the absorbed left subtrajectory is the connecting state. The second
demands that the multiplier of continuity is the gradient of the right value at the connecting state. Together they
give the optimal parametric costate and connecting state
\begin{equation}
\begin{aligned}
\lambda^*_{\tau,t\to t'}(x_t,\lambda_{t'})&=L_{\tau,t\to t',x}x_t+L_{\tau,t\to t',\lambda}\lambda_{t'}+l_{\tau,t\to t',0},\label{eq:lambdamap}\\
x^*_{\tau,t\to t'}(x_t,\lambda_{t'})&=X_{\tau,t\to t',x}x_t+X_{\tau,t\to t',\lambda}\lambda_{t'}+x_{\tau,t\to t',0},
\end{aligned}
\end{equation}
and the combined tuple. 

\begin{theorem}\label{thm:comb}
Absorption followed by the junction combines two consecutive LQ costate-conditioned value tuples into the
costate-conditioned value tuple of $[t,t']$, and a costate-conditioned and a conventional value tuple into the
conventional value tuple at $t$. At the optimum, with $\hat\nu_\tau$ supplied by the
parent, the multiplier of the right domain is
\begin{equation}\label{eq:nucompfull}
\nu^*_{\tau,t\to t'}=\tilde T_\tau\tilde\nu^*_{\tau,t\to t'}(x_t,\lambda^*_\tau)+\hat T_\tau\hat\nu_\tau,
\end{equation}

The effective costate, $\mu^*_\tau=\lambda^*_{\tau,t\to t'}+\Wc_{\tau\to t',x}^\T\nu^*_{\tau,t\to t'}$, equals the costate, $\pi_\tau$, of
\eqref{eq:lagrangian}.
\end{theorem}

\subsubsection{Interpretation of the multipliers} Three multipliers are associated with the cut $\tau$. Their interpretation is essential to understand the principle bookkeeping of the method. 
\begin{itemize}[leftmargin=*]
\item The continuity multiplier, $\lambda^*_{\tau,t\to t'}$, corresponds with the gradient of the downstream value
$\Vc_{\tau\to t'}$ at the connecting state, i.e.\ the price of $x_\tau$, \emph{without} the downstream domain.
\item The multiplier, $\nu^*_{\tau,t\to t'}$, corresponds with the downstream domain $\wc_{\tau\to t'}$. The steerable part is determined at this node by the
absorption step, because the upstream subtrajectory can satisfy these rows. The unsteerable part, $\hat{\nu}_\tau$, is determined by an ancestor,
because only the decision that determines $x_t$ can settle them. This is equivalent to the decomposition $\eta=\tilde T\tilde\eta+\hat T\hat\eta$ of
Proposition~\ref{prop:pqp}.
\item Th effective costate, $\mu^*_\tau$, is the price of $x_\tau$ \emph{including} the downstream domain,
\begin{equation}\label{eq:musplit}
\mu^*_\tau=\lambda^*_\tau+\tilde W_{\tau\to t',x}^\T\tilde\nu^*_\tau+\hat W_{\tau\to t',x}^\T\hat\nu_\tau .
\end{equation}
Seen from the upstream subtrajectory,
$\mu^*_\tau$ is the slope of its terminal cost, hence its terminal costate. Seen from the downstream subtrajectory, $x_\tau$ is the parameter and by \eqref{eq:pqpsens}, the price of $x_\tau$ is
$\mu^*_\tau=\nabla_x\Vc_{\tau\to t'}+\Wc_{\tau\to t',x}^\T\nu^*_{\tau}$. 
\end{itemize}

The combination procedure discussed above is again associative and supports parallel-in-time solution of the problem. We discuss the procedure in analogy with the parallel treatment from section \ref{sec:parsolstatecond} focussing on the primal solution. Recovery of the dual solution will be detailled in another section.



\subsection{Parallel solution}\label{sec:costatepar}

Assume the horizon of \eqref{eq:ocp} is again cut at time instants $t_0=0<t_1<\dots<t_{N-1}<t_N=T$. 
\begin{itemize}
    \item In an \emph{initial pass} the conditional, $\cV_{t_i\to t_{i+1}}$, $i\in\N_{0:N-2}$, and conventional, $\mathcal{V}_{t_{N-1}}$, tuples are computed independently on different cores.
    \item In an \emph{upward sweep} every
node combines its children by Theorem~\ref{thm:comb} and stores the corresponding maps. 
\item In the \emph{downward sweep} a node with junction $\tau$ receives the optimal initial state, $x^*_t$, the
effective costate, $\mu^*_{t'}$, and the multiplier
$\nu^*_{t}$ of its (absorbed) domain. Since $\wc_{t\to t'}$ stacks the domain $\wc_{t\to\tau}$ of the upstream child and the transferred rows
\eqref{eq:newrows}, the multiplier splits accordingly,
\begin{equation}\label{eq:nustack}
\nu^*_{t}=\begin{pmatrix}\nu^*_{t\to\tau}\\\hat\nu_\tau\end{pmatrix}.
\end{equation}
Then the node evaluates the variables at the cut.  First the node evaluates $x^*_\tau$ and $\lambda^*_\tau$ by the maps
\eqref{eq:lambdamap} with $\lambda_{t'}\leftarrow \mu_{t'}^*$. Second the node evaluates \eqref{eq:numap} and \eqref{eq:mueff}
\begin{equation}
\begin{aligned}
\tilde\nu^*_\tau={}&N_{\tau,t\to t',x}x_t^*+N_{\tau,t\to t',\lambda}\lambda_{\tau}^*+n_{\tau,t\to t',0},\label{eq:nutdown}\\
\nu^*_{\tau}={}&\tilde T_\tau\tilde\nu^*_\tau+\hat T_\tau\hat\nu_\tau,\\
\mu^*_\tau={}&\lambda_\tau^* +\Wc_{\tau\to t',x}^\T\nu_{\tau}^*. 
\end{aligned}
\end{equation}

Finally, the node then passes $(x^*_t,\mu^*_\tau,\nu^*_{t})$ upstream and
$(x^*_\tau,\mu^*_{t'},\nu^*_{\tau})$ downstream.
\item A \emph{final pass} simulates every
subtrajectory in closed loop with 
\begin{equation}\label{eq:lqfinal}
u^*_t=U_{t\to t_{i+1},x}x^*_t+U_{t\to t_{i+1},\lambda}\mu^*_{t_{i+1}}+u_{t\to t_{i+1},0}.
\end{equation}
and $t\in[t_i,t_{i+1}]$ producing the primal solution, consisting again of $x_{0:T}^*$ and $u_{0:T-1}^*$.
\end{itemize}

%


\subsection{Dual solution recovery}\label{sec:dual}

The procedure discussed above so far only recovers the primal solution. Now we discuss how to retrieve its dual, consisting of the multipliers $\pi_{1:T}^*$ and $\gamma_{0:T}^*$. The full procedure, including the dual recovery, is summarized in Algorithm \ref{alg:main}.

To recover the dual variables, we rely on Proposition~\ref{prop:pqp} which suggests that every constraint is {used} at exactly one elimination step, which then determines the associated multiplier. 
At stage $t$ of leaf $i$, with $t'=t_{i+1}$, the p-QP deals with constraint rows, $g_t$, and the domain rows,
$\wc_{t+1\to t'}(x_{t+1})$. The corresponding multiplier decomposes as $\eta_t=(\gamma_t,\nu_{t+1})$. The rows that are transferred upstream are the domain rows at $t$, with multiplier $\nu_{t}$. By definition, multipliers $\pi_t$ determine the price of the parameter $x_t$ and can be retrieved through the gradient of the value tuple.

\begin{proposition}\label{prop:mult}
Along the optimal trajectory, with $T_t$ the row transformation of the stage p-QP and $\tilde\eta_t$ the multiplier
\eqref{eq:pqpmult2} of the rows is uses,
\begin{equation}
\begin{aligned}
\eta_t&=T_t\begin{pmatrix}\tilde\eta_t\\\nu_{t}\end{pmatrix},\label{eq:mr1}\\
\pi_{t+1}&=\nabla_x\Vc_{t+1\to t'}(x_{t+1},\mu_{t'}^*)+\Wc_{t+1\to t',x}^\T\nu_{t+1}.
\end{aligned}
\end{equation}
\end{proposition}

In the downward sweep, at the cuts, $\nu^*_{\tau}$ corresponds with the multiplier of the domain of the downstream child, and $\mu^*_\tau$ corresponds with the costate $\pi_\tau$.

For all stages of a leaf, given $x^*_{t_i}$, $\mu^*_{t_{i+1}}$ and $\nu^*_{t_i}$, one forward sweep produces $u_t^*$, $x_{t+1}^*$, $\gamma_t^*$, $\nu_{t+1^*}$ and $\pi_{t+1}^*$. At the end of the last or final leaf we have that $\nu_T=\gamma_T$. 



\begin{algorithm}[t]
\caption{Costate-conditioned parallel solution of \eqref{eq:ocp}}\label{alg:main}
\begin{algorithmic}[1]
\Require data \eqref{eq:lqdata}, $x_0$, cuts $0=t_0<\dots<t_{N-1}<T$
\State \textbf{parfor} leaves: $\cV_{t\to t_{i+1}}$, $u^*_{t\to t_{i+1}}$, $t\in\N_{t_i:t_{i+1}-1}$ (last leaf is conventional: $\mathcal{V}_{t_{N-1}}$, $u^*_{t_{N-1}}$), pass $\cV_{t_i\to t_{i+1}}$ to the tree
\State \textbf{for} each level, bottom-up, \textbf{parfor} each node: absorption, junction, store maps $N$, $L$, $X$
\State \textbf{for} each level, top-down, \textbf{parfor} each node: $x^*_\tau$, $\lambda^*_\tau$, $\nu^*_{\tau\to t'}$, $\mu^*_\tau$, pass
$(x^*_t,\mu^*_\tau,\nu^*_{t})$ left and $(x^*_\tau,\mu^*_{t'},\nu^*_{\tau})$ right
\State \textbf{parfor} leaves: forward sweep \eqref{eq:lqfinal}, \eqref{eq:mr1}, recovering the full primal-dual solution $\{x_{1:T}^*,u_{0:T-1}^*,\pi_{1:T}^*,\gamma_{1:T}^*\}$
\end{algorithmic}
\end{algorithm}




\section{Complexity}\label{sec:complexity}

In conclusion we provide an analysis of the complexity of the algorithm. One may note in advance that the goal of costate-conditioning was never to reduce the total operation count w.r.t. to state-conditionining but to improve the numerical behaviour, as will be demonstrated in Section~\ref{sec:numerics}.


Let $\omega$ be the work of one stage of the sequential recursion, $\rho\,\omega$ that of one stage of a conditioned recursion in a leaf, and
$\kappa\,\omega$ that of one combination. With $T$ stages and $N$ subtrajectories of equal length the bound on the speed-up of the linear
solve is
\begin{equation}\label{eq:speedup}
\varphi(N)=\frac{W_\text{seq}}{W_\text{par}}\approx\frac{T\omega}{\frac{T}{N}\rho \omega+\lceil\log_2N\rceil\kappa\omega} \approx \frac{N}{\rho}
\end{equation}
The second approximation is valid when the first term in the denominator dominates of the second term, so generally for very large $T$. 

By Proposition~\ref{prop:leaf}, $\rho > 1$ is a small constant independent of $T$. Note that this model does not account for synchronisation and data movement.


\section{Numerical results}\label{sec:numerics}

All methods were implemented in MATLAB in double precision with a common p-QP routine. Rank decisions use a relative singular-value threshold of $10^{-10}$. The reference is a sparse direct primal-dual solve of the full KKT system. We compare the costate-conditioned tree (Algorithm~\ref{alg:main}), the state-conditioned tree \cite{lefebvre2026temporal} and a sequential recursion (i.e. Algorithm~\ref{alg:main} with a single leaf). 

\subsection{Test problems}

\subsubsection{Random problems}\label{sec:numverif}

We first evaluate the algorithm on random problems with $n_x=4$, $n_u=2$, $T=16$, cuts at $t\in\{0,4,8,12\}$, mixed state-input constraints at $t\in\{2,5,9,13\}$, a
rank-2 state-only constraint at $t=8$, a fully pinned state at $t=12$ and a rank-3 terminal constraint. These settings exercise every
mechanism in the proposed algorithm. For completeness, besides a balanced recombination (which is most cost effective), we also evaluate left and right fold recombinations.


\begin{table}[t]\centering\footnotesize
\setlength{\tabcolsep}{3pt}
\caption{Largest absolute deviation of the primal-dual solution obtained with various algorithms from the direct KKT solution. A distinction is made for problems where the stage Hessian was convex or indefinite.}\label{tab:verif}
\begin{tabular}{@{}lcc@{}}
\toprule
 & convex & indefinite \\
\midrule
sequential recursion & $2.0\cdot10^{-13}$ & $3.3\cdot10^{-12}$ \\
costate tree, balanced & $4.0\cdot10^{-13}$ & $3.3\cdot10^{-13}$ \\
costate tree, left fold & $7.4\cdot10^{-13}$ & $1.3\cdot10^{-12}$ \\
costate tree, right fold & $4.0\cdot10^{-13}$ & $4.2\cdot10^{-13}$ \\
state-conditioned tree & $4.0\cdot10^{-12}$ & $3.6\cdot10^{-13}$ \\
\bottomrule
\end{tabular}
\end{table}

\subsubsection{Waypoint tracking}\label{sec:numtrack}

In a second experiment we recycle the tracking problem of \cite{lefebvre2026temporal}. The system is determined by 3D Cartesian coordinates and velocities $x_t=(p_t,v_t)\in\R^6$ of a point mass, a thrust input $u_t\in\R^3$ and double-integrator
dynamics with step $\Delta t=0.1$. We impose quadratic costs on the velocity and input and introduce $K+1$ position waypoint constraints from a Brownian walk with unit steps and
$L$ stages between every two waypoints bringing the total length to $T=KL$. The terminal state is pinned.
Cuts are placed at the waypoints.

\subsubsection{Chain of masses}\label{sec:conditioning}

For the tracking problem, which has controlability index $2$, all trees are as accurate as the sequential recursion. To illustrate the reachability problem we introduce a final experiment determined by a chain of masses driven by a
force exerted on the last mass which has controllability index $n_x$ with $n_x$ twice the number of masses. The total horizon is $T=256$ and the terminal state is pinned. We solve the problem for varying subtrajectory lengths and varying state size.

\begin{figure}[h!]\centering
\includegraphics[width=.95\columnwidth]{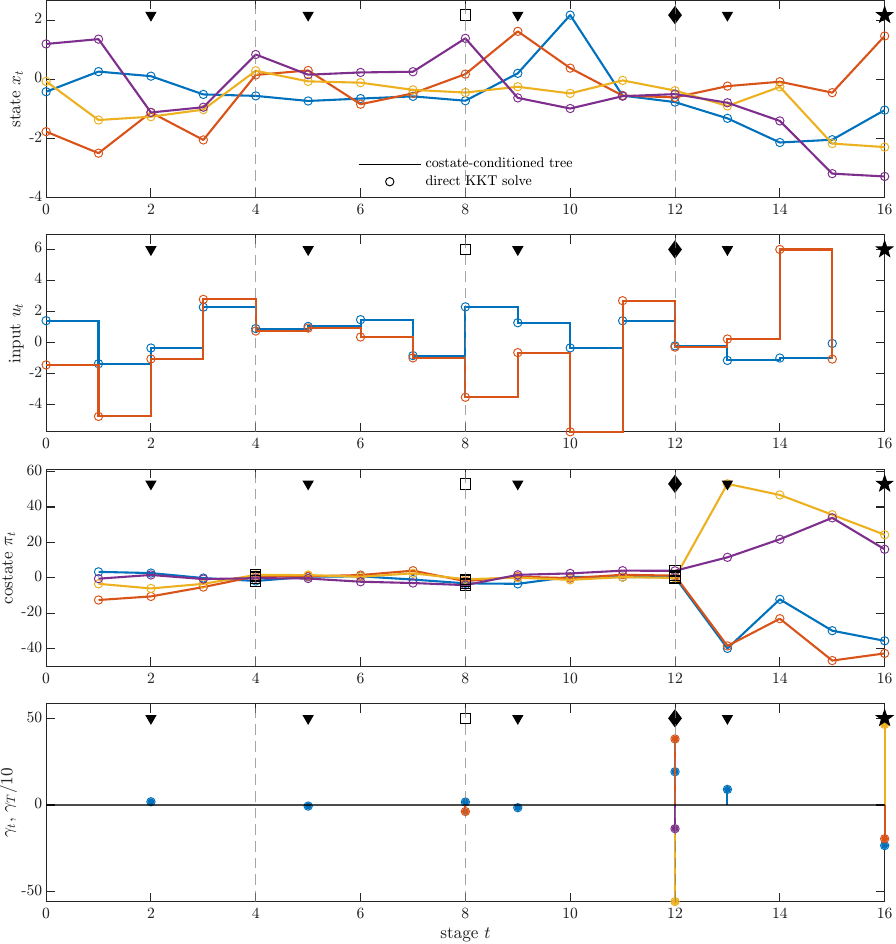}
\caption{Illustration of the solution of a random test problem. First two panels show the primal solution, last two panels the dual. The markers on top reveal the presence of constraints, i.e. mixed ($\blacktriangledown$), rank-2 state-only ($\square$),
fully pinned state ($\blacklozenge$) and terminal ($\bigstar$). The full lines, circles and dashed lines represent respectively the costate-conditioned solution, the direct KKT solve and the location of the cuts. The black squares in the third panel indicate the effective
costates, $\mu^*_{t_i}$, generated by the downward sweep. }
\label{fig:random}
\vspace{10pt}
\includegraphics[width=.95\columnwidth]{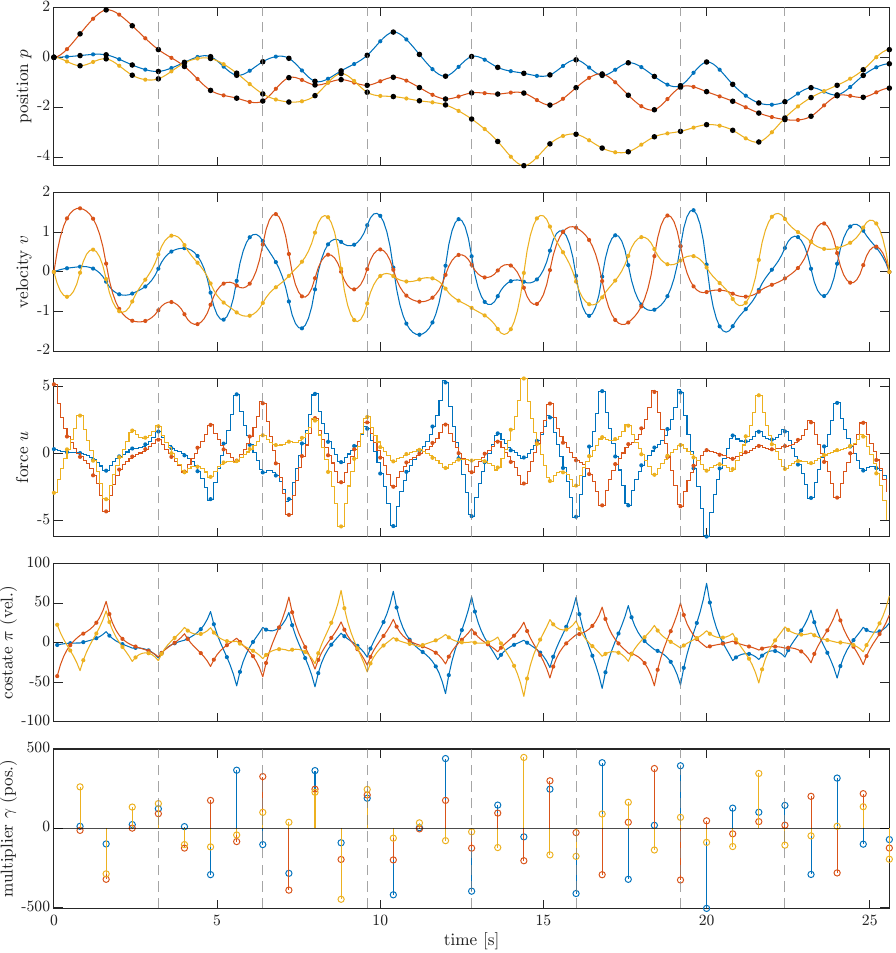}
\caption{Illustration of the solution of the tracking problem with $K=32$, $L=8$, $N=8$. Interpretation of the figure is similar to those in Fig. \ref{fig:random}.}
\label{fig:trackprof}
\end{figure}

\begin{figure}[t]\centering
\includegraphics[width=.95\columnwidth]{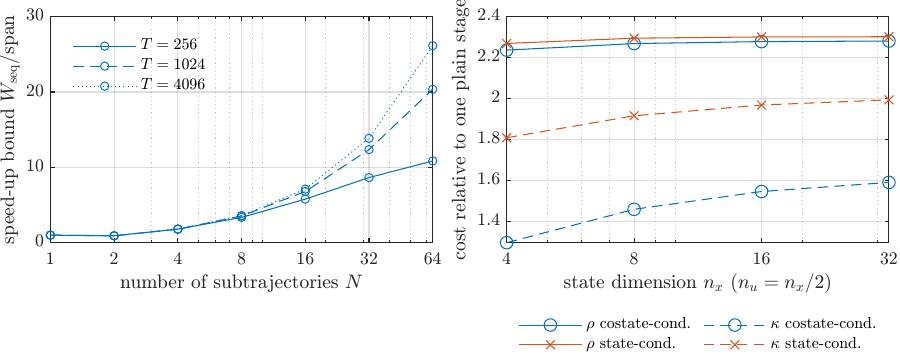}
\caption{Left: bound on the speed-up of the linear solve with the costate-conditioned tree for the tracking problem
($K=64$, cuts at waypoints). Right:
measured $\rho$ and $\kappa$ against the state dimension on random problems.}
\label{fig:complexity}
\vspace{10pt}

\includegraphics[width=.95\columnwidth]{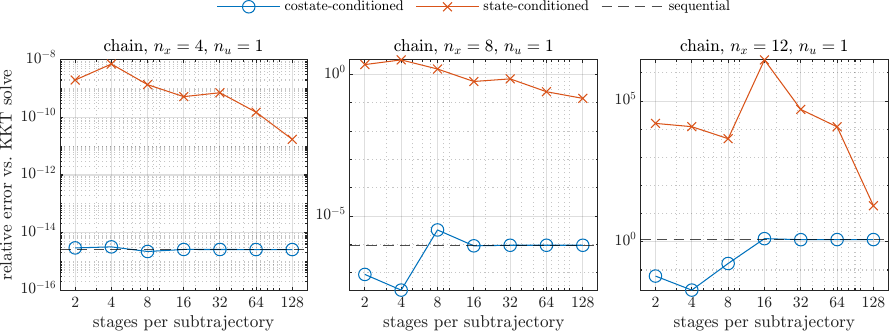}
\caption{Chain of masses ($n_u=1$, controllability index $n_x$): relative error of states and inputs against the direct
KKT solve versus the number of stages per subtrajectory. Dashed: sequential recursion.}
\label{fig:chain}
\vspace{-10pt}
\end{figure}

\subsection{Results}

Table~\ref{tab:verif}
lists the largest deviation from the reference for the random problem. All methods agree with the reference to within $4\cdot10^{-12}$, and the costate-conditioned tree to within $1.3\cdot10^{-12}$ in all quantities, including all multipliers. Fig. \ref{fig:random} visualizes the solution of one of the random problems.

Fig.~\ref{fig:trackprof} illustrates the solution for a random set of waypoints with $K=32$, $L=8$ and $N=8$. The full primal-dual solution agrees with the
reference to a relative $2\cdot10^{-13}$. 

Second we visualize the theoretical speed-up for the tracking problem. The result is visualized in Fig.~\ref{fig:complexity}, left panel. For $N=64$ the speed-up bound is $26$ with $T=4096$, $20$ for
$T=1024$ and $11$ for $T=256$. 

The factors, $\rho$, and, $\kappa$, were evaluated empirically for varying state dimension on variations of the random problem. The results are visualized in Fig.~\ref{fig:complexity}, right panel. The factor, $\rho$, is $\rho\approx2.25$, irrespective of the state dimension and the conditioning type. The factor, $\kappa$, ranges from $1.3$--$2.2$ with a constant difference of $0.4$ between state- and costate conditioned combinations.

\section{Conclusion}\label{sec:conclusion}

This work is tailored to parallel-in-time primal-dual solution of equality constrained linear-quadratic optimal control problems. This problem is recognized as the computational bottleneck in gradient based trajectory optimization algorithms. 

Dual-regularised parallel-in-time methods \cite{jallet2024parallel,sousapinto2025dual} eliminate all stagewise equality constraints before
the parallel step and remove all rank decisions, but each Newton step satisfies the linearised constraints only up to the regularisation. Others methods \cite{lefebvre2026temporal,laine2019parallelizing,vanroye2023fatrop} aim to solve the problem exactly. The existing approaches relied on state-conditioned value tuples which, although elegant, results in numerical disadvantages. 

This work proposed to condition the value tuples on the costate. This results in a free-end constrained Riccati recursions without reachability rows. The latter property makes the proposed method numerically advantageous over state-conditioned tuples and offers additional benefits when recovering the dual solution.

The complete primal--dual solution follows at
small additional cost, which makes the method a candidate linear solver for IPOPT-type trajectory optimisation. Future
work concerns implementations on multi-core and GPU hardware and integration in a trajectory optimization solver.


\section*{AI declaration}
Large Language Models (LLMs) were used as a general-purpose research assistant, under the direction of the author, for drafting and editing text, for implementing and running experiments, for literature search, and as a critical reviewer of derivations and experimental protocols. All theoretical claims, experimental designs, and conclusions were reviewed and verified by the authors, who take full responsibility for the content of the paper.

\small
\balance
\bibliographystyle{IEEEtran}
\bibliography{references}

\onecolumn
\clearpage
\appendices

\section{Proofs}\label{app:proofs}

\begin{proofof}{Proposition~\ref{prop:pqp}}
(i) With \eqref{eq:pqpmult1} and \eqref{eq:pqpT}, $B_u^\T\eta=(\tilde T^\T B_u)^\T\tilde\eta+(\hat T^\T B_u)^\T\hat\eta=\tilde B_u^\T\tilde\eta$, so
stationarity in $u$ is \eqref{eq:pqpmult2} and does not involve $\hat\eta$. Since $u^*(x)$ minimises $a(x,\cdot)$ on the
affine set of inputs that satisfy the used rows, $\nabla_ua(x,u^*)$ is orthogonal to the null space of $\tilde B_u$, hence
in the range of $\tilde B_u^\T$; as $\tilde B_u$ has full row rank, $\tilde\eta$ exists and is unique. (ii) The policy satisfies the
used rows for every $x$, so $\tilde B_uU_x=-\tilde B_x$ with $\tilde B_x=\tilde T^\T B_x$. By the chain rule and (i),
$\nabla_x\hat a=\nabla_xa+U_x^\T\nabla_ua=\nabla_xa-U_x^\T\tilde B_u^\T\tilde\eta=\nabla_xa+\tilde B_x^\T\tilde\eta$, and
$B_x^\T\eta=\tilde B_x^\T\tilde\eta+\hat B_x^\T\hat\eta$. Combining both gives \eqref{eq:pqpsens}.
\end{proofof}

\begin{lemma}\label{lem:conj}
For $x_t\in\cX_{t\to t'}$ and $(x_t,x_{t'})\in\mathcal{X}_{t\to t'}$, respectively,
\begin{equation}
\begin{aligned}
\Vc_{t\to t'}(x_t,\lambda_{t'})={}&\min_{x_{t'}}\ V_{t\to t'}(x_t,x_{t'})+\lambda_{t'}^\T x_{t'}\notag\\
&\ \text{s.t. }0=w_{t\to t'}(x_t,x_{t'}),\label{eq:conj1}\\
V_{t\to t'}(x_t,x_{t'})={}&\ext_{\lambda_{t'}}\ \Vc_{t\to t'}(x_t,\lambda_{t'})-\lambda_{t'}^\T x_{t'}.
\end{aligned}
\end{equation}
\end{lemma}

\begin{proof}
Splitting the minimisation in \eqref{eq:defcostate} into one over $x_{t'}$ and one over the inputs with $x_{t'}$ fixed gives
\eqref{eq:conj1}, since the inner problem is \eqref{eq:cond}. For \eqref{eq:conj1}, the right-hand side is the stationary value
of the Lagrangian of \eqref{eq:cond} in which only the constraint is dualised with multiplier $\lambda_{t'}$; its stationarity
conditions together with the remaining constraints are the optimality conditions of \eqref{eq:cond}. On
$\mathcal{X}_{t\to t'}$ the stationary $\lambda_{t'}$ is unique up to components orthogonal to the reachable terminal states, which do
not change the value.
\end{proof}

\begin{proofof}{Theorem~\ref{thm:leaf}}
Lemma~\ref{lem:bellman} on $[t,t']$ with terminal cost $\lambda_{t'}^\T x_{t'}$, no terminal constraint and $\lambda_{t'}$ as a parameter.
\end{proofof}

\begin{proofof}{Proposition~\ref{prop:leaf}}
(i) By induction: if $\Wc_{t+1\to t',\lambda}=0$, the constraint element \eqref{eq:cp4} does not depend on $\lambda_{t'}$, so the
transferred rows do not either and $U_{\perp,\lambda}=0$. (ii) The costate enters the p-QP data only through the blocks
\eqref{eq:cq8}--\eqref{eq:cq10}, which do not affect $B_u$, $B_x$, $A_{uu}$, $A_{ux}$, $A_{xx}$; at $\lambda_{t'}=0$ the problem is the
free-end problem. (iii) With $U_{\perp,\lambda}=0$, $U_{t\to t',\lambda}=-T_z\bar A_{zz}^{-1}\bar A_{z\lambda}$ with
$\bar A_{z\lambda}=T_z^\T\Qc_{t\to t',\lambda u}^\T$, and \eqref{eq:pqpval1} gives
$\Vc_{t\to t',\lambda\lambda}=\Vc_{t+1\to t',\lambda\lambda}-\bar A_{z\lambda}^\T\bar A_{zz}^{-1}\bar A_{z\lambda}$.
\end{proofof}

\begin{proofof}{Theorem~\ref{thm:comb}}
Dualise continuity with $\lambda_\tau$ and the right domain, imposed on $x_\tau$, with $\nu_{\tau\to t'}$. The Lagrangian separates
into a left part, the subtrajectory $[t,\tau]$ with terminal cost $\mu_\tau^\T x_\tau+\nu_{\tau\to t'}^\T\wc_{\tau\to t',0}$, and a right
part, the subtrajectory $[\tau,t']$ with the additional initial term $-\lambda_\tau^\T x_\tau$. Their stationary values are
$\Vc_{t\to\tau}(x_t,\mu_\tau)+\nu_{\tau\to t'}^\T\wc_{\tau\to t',0}$ and $-\lambda_\tau^\T x_\tau+\Vc_{\tau\to t'}(x_\tau,\lambda_{t'})$, and as in
Lemma~\ref{lem:conj} the stationary value over $(x_\tau,\lambda_\tau,\nu_{\tau\to t'})$ is the value of the coupled problem. By
\eqref{eq:conj1} and Theorem~\ref{thm:lvcomb} this is the costate-conditioned value of $[t,t']$, respectively the value at $t$.

\emph{Absorption.} By \eqref{eq:termstate} the left terminal state is
$x_\tau=\Vc_{t\to\tau,\lambda x}x_t+\vc_{t\to\tau,\lambda}+\Vc_{t\to\tau,\lambda\lambda}\mu_\tau$, and stationarity in $\nu_{\tau\to t'}$ is
$\Wc_{\tau\to t',x}x_\tau+\wc_{\tau\to t',0}=0$. Along $\tilde T_\tau$ it is a linear equation for $\tilde\nu_\tau$ with matrix
$S_\tau=\tilde W_{\tau\to t',x}\Vc_{t\to\tau,\lambda\lambda}\tilde W^\T_{\tau\to t',x}$, which is nonsingular by the choice of $\tilde T_\tau$ if
$\Vc_{t\to\tau,\lambda\lambda}\preceq0$ (Proposition~\ref{prop:leaf}(iii)) and is assumed nonsingular otherwise; $\hat\nu_\tau$ does not enter because
$\Vc_{t\to\tau,\lambda\lambda}\hat W^\T_{\tau\to t',x}=0$. Along $\hat T_\tau$ it is \eqref{eq:newrows}, which does not involve $\nu_{\tau\to t'}$; the
Lagrangian is linear in $\hat\nu_\tau$ with this residual as coefficient, so \eqref{eq:newrows} becomes a domain row and
$\hat\nu_\tau$ its multiplier, determined where $x_t$ is determined.

\emph{Junction.} Stationarity in $(\lambda_\tau,x_\tau)$ of the remaining problem gives \eqref{eq:stat1}.

\emph{Associativity.} The combined tuple is the costate-conditioned tuple of $[t,t']$, which by Lemma~\ref{lem:conj} is a
different representation of the state-conditioned tuple; the combination of the latter is associative
\cite[Corollary~4]{lefebvre2026temporal}.

\emph{Costate.} The left subtrajectory with terminal cost $\mu^{*\T}_\tau x_\tau$ has the same stationarity conditions
as the full problem restricted to $[t,\tau]$, with terminal costate $\mu^*_\tau$; hence $\mu^*_\tau=\pi_\tau$.
\end{proofof}

\begin{corollary}\label{cor:policy}
For $t<\tau\le T$, $u^*_t(x_t)=u^*_{t\to\tau}(x_t,\mu^*_{\tau,t}(x_t))$, with $\mu^*_{\tau,t}$ the effective costate of the
combination with the conventional tuple $\mathcal{V}_\tau$. Then the term $\hat W^\T\hat\nu_\tau$ does not affect the inputs and may be omitted when only the primal solution is required.
\end{corollary}
\begin{proof}
By Theorem~\ref{thm:comb} the optimal inputs on $[t,\tau)$ are those of the left policy at the effective costate. The rows
associated with $\hat\nu_\tau$ fix a function of $x_\tau$ that the inputs cannot move, so the corresponding term of the terminal
cost does not depend on the inputs.
\end{proof}

\begin{proofof}{Proposition~\ref{prop:mult}}
Equation \eqref{eq:mr1} is \eqref{eq:pqpmult1}--\eqref{eq:pqpmult2} with $\hat\eta_t=\nu_{t\to t'}$. Expanding the $Q$-tuple, the
stationarity condition \eqref{eq:pqpmult2} of the p-QP in $u_t$ is the $u_t$-row of the optimality conditions of
\eqref{eq:ocp} with $\pi_{t+1}$ given by \eqref{eq:mr1}. The left-hand side of \eqref{eq:pqpsens} for the p-QP at stage $t+1$ is
$\nabla_xc_{t+1}+G_{t+1,x}^\T\gamma_{t+1}+F_{t+1,x}^\T\pi_{t+2}$ and its right-hand side is \eqref{eq:mr1}, so the $x_{t+1}$-row
holds; at $t+1=t'$, \eqref{eq:mr1} gives $\pi_{t'}=\mu$. In the tree, the rows along $\tilde T_\tau$ are used in the absorption and
have the multiplier \eqref{eq:numap} at $\lambda^*_\tau$, those along $\hat T_\tau$ are passed to the parent and receive $\hat\nu_\tau$,
and multipliers transform back with $T_\tau$, which gives in $\nu_\tau$ in \eqref{eq:nutdown}. By Theorem~\ref{thm:comb}, $\mu_\tau$ is
$\pi_\tau$.
\end{proofof}

\section{Null-space p-QP}\label{app:pqp}

The p-QP \eqref{eq:pqp} has the cost
\begin{equation}
a(x,u)=\tfrac12\begin{pmatrix}1\\x\\u\end{pmatrix}^{\!\T}
\begin{pmatrix}2a_0&\bul&\bul\\a_x&A_{xx}&\bul\\a_u&A_{ux}&A_{uu}\end{pmatrix}
\begin{pmatrix}1\\x\\u\end{pmatrix}.
\end{equation}
With $T$ as in \eqref{eq:pqpT}, let
\begin{align}
\tilde B_x&=\tilde T^\T B_x,\\
\tilde b_0&=\tilde T^\T b_0,
\end{align}
and let $[T_y\ T_z]$ be orthogonal with $\tilde B_uT_z=0$. Then
\begin{align}
U_{\perp,x}&=-T_y(\tilde B_uT_y)^{-1}\tilde B_x,\\
u_{\perp,0}&=-T_y(\tilde B_uT_y)^{-1}\tilde b_0,\\
\bar A_{zz}&=T_z^\T A_{uu}T_z,\\
\bar A_{zx}&=T_z^\T(A_{ux}+A_{uu}U_{\perp,x}),\\
\bar a_z&=T_z^\T(a_u+A_{uu}u_{\perp,0}),\\
U_x&=U_{\perp,x}-T_z\bar A_{zz}^{-1}\bar A_{zx},\\
u_0&=u_{\perp,0}-T_z\bar A_{zz}^{-1}\bar a_z,\\
\hat A_{xx}&=A_{xx}+A_{ux}^\T U_x+U_x^\T A_{ux}+U_x^\T A_{uu}U_x,\label{eq:pqpval1}\\
\hat a_x&=a_x+A_{ux}^\T u_0+U_x^\T(a_u+A_{uu}u_0),\label{eq:pqpval2}
\end{align}
where $\hat A_{xx}$ and $\hat a_x$ are the coefficients of $\hat a(x)$ and $\bar A_{zz}$ is the pivot. The multiplier of the used rows is
\begin{equation}\label{eq:pqpmult3}
\tilde\eta=-(\tilde B_uT_y)^{-\T}T_y^\T\big(a_u+A_{ux}x+A_{uu}u^*\big).
\end{equation}

\section{LQ coefficients of the state-conditioned tuples}\label{app:lvlq}

The $Q$-tuple of Theorem~\ref{thm:condrec} has the cost coefficients
\begin{align}
q_{t\to t',x}&=c_{t,x}+F_{t,x}^\T\big(v_{t+1\to t',x}+V_{t+1\to t',xx}f_{t,0}\big),\\
q_{t\to t',u}&=c_{t,u}+F_{t,u}^\T\big(v_{t+1\to t',x}+V_{t+1\to t',xx}f_{t,0}\big),\\
q_{t\to t',x'}&=v_{t+1\to t',x'}+V_{t+1\to t',x'x}f_{t,0},\\
Q_{t\to t',xx}&=C_{t,xx}+F_{t,x}^\T V_{t+1\to t',xx}F_{t,x},\\
Q_{t\to t',ux}&=C_{t,ux}+F_{t,u}^\T V_{t+1\to t',xx}F_{t,x},\\
Q_{t\to t',uu}&=C_{t,uu}+F_{t,u}^\T V_{t+1\to t',xx}F_{t,u},\\
Q_{t\to t',x'x}&=V_{t+1\to t',x'x}F_{t,x},\\
Q_{t\to t',x'u}&=V_{t+1\to t',x'x}F_{t,u},\\
Q_{t\to t',x'x'}&=V_{t+1\to t',x'x'},
\end{align}
and the constraint coefficients
\begin{align}
P_{t\to t',x}&=\begin{pmatrix}G_{t,x}\\W_{t+1\to t',x}F_{t,x}\end{pmatrix},\\
P_{t\to t',u}&=\begin{pmatrix}G_{t,u}\\W_{t+1\to t',x}F_{t,u}\end{pmatrix},\\
P_{t\to t',x'}&=\begin{pmatrix}0\\W_{t+1\to t',x'}\end{pmatrix},\\
p_{t\to t',0}&=\begin{pmatrix}g_{t,0}\\W_{t+1\to t',x}f_{t,0}+w_{t+1\to t',0}\end{pmatrix}.
\end{align}
The p-QP with variable $u_t$, parameter $(x_t,x_{t'})$, this cost and this constraint yields \eqref{eq:lvw}--\eqref{eq:lvpol}. The
boundary condition \eqref{eq:lvinit1} sets all cost coefficients to zero and
\begin{align}
W_{t'\to t',x}&=I,\\
W_{t'\to t',x'}&=-I,\\
w_{t'\to t',0}&=0.
\end{align}
The conventional recursion omits $x_{t'}$ and starts from
\begin{align}
v_{T,x}&=c_{T,x},\\
V_{T,xx}&=C_{T,xx},\\
W_{T,x}&=G_{T,x},\\
w_{T,0}&=g_{T,0}.
\end{align}
The combination \eqref{eq:lvc2} is the p-QP with variable $x_\tau$, parameter $(x_t,x_{t'})$ and
\begin{align}
A_{uu}&=V_{t\to\tau,x'x'}+V_{\tau\to t',xx},\\
A_{ux}&=\begin{pmatrix}V_{t\to\tau,x'x}&V_{\tau\to t',x'x}^\T\end{pmatrix},\\
a_u&=v_{t\to\tau,x'}+v_{\tau\to t',x},\\
B_u&=\begin{pmatrix}W_{t\to\tau,x'}\\W_{\tau\to t',x}\end{pmatrix},\\
B_x&=\begin{pmatrix}W_{t\to\tau,x}&0\\0&W_{\tau\to t',x'}\end{pmatrix},\\
b_0&=\begin{pmatrix}w_{t\to\tau,0}\\w_{\tau\to t',0}\end{pmatrix},
\end{align}
whose policy is \eqref{eq:lvX}.

\section{LQ coefficients of the costate-conditioned tuples}\label{app:costatelq}

\emph{Leaf recursion.} The $Q$-tuple of Theorem~\ref{thm:leaf} has the cost coefficients
\begin{align}
\qc_{t\to t',x}&=c_{t,x}+F_{t,x}^\T\big(\vc_{t+1\to t',x}+\Vc_{t+1\to t',xx}f_{t,0}\big),\\
\qc_{t\to t',u}&=c_{t,u}+F_{t,u}^\T\big(\vc_{t+1\to t',x}+\Vc_{t+1\to t',xx}f_{t,0}\big),\\
\qc_{t\to t',\lambda}&=\vc_{t+1\to t',\lambda}+\Vc_{t+1\to t',\lambda x}f_{t,0},\\
\Qc_{t\to t',xx}&=C_{t,xx}+F_{t,x}^\T\Vc_{t+1\to t',xx}F_{t,x},\\
\Qc_{t\to t',ux}&=C_{t,ux}+F_{t,u}^\T\Vc_{t+1\to t',xx}F_{t,x},\\
\Qc_{t\to t',uu}&=C_{t,uu}+F_{t,u}^\T\Vc_{t+1\to t',xx}F_{t,u},\\
\Qc_{t\to t',\lambda x}&=\Vc_{t+1\to t',\lambda x}F_{t,x},\label{eq:cq8}\\
\Qc_{t\to t',\lambda u}&=\Vc_{t+1\to t',\lambda x}F_{t,u},\\
\Qc_{t\to t',\lambda\lambda}&=\Vc_{t+1\to t',\lambda\lambda},\label{eq:cq10}
\end{align}
and the constraint coefficients
\begin{align}
\Pc_{t\to t',x}&=\begin{pmatrix}G_{t,x}\\\Wc_{t+1\to t',x}F_{t,x}\end{pmatrix},\\
\Pc_{t\to t',u}&=\begin{pmatrix}G_{t,u}\\\Wc_{t+1\to t',x}F_{t,u}\end{pmatrix},\\
\pc_{t\to t',0}&=\begin{pmatrix}g_{t,0}\\\Wc_{t+1\to t',x}f_{t,0}+\wc_{t+1\to t',0}\end{pmatrix}.\label{eq:cp4}
\end{align}
For the forward multiplier recursion, \eqref{eq:pqpmult3} for the stage p-QP reads
\begin{align}
\tilde\eta_t={}&-(\tilde B_uT_y)^{-\T}T_y^\T\big(\qc_{t\to t',u}+\Qc_{t\to t',ux}x_t\notag\\
&+\Qc_{t\to t',\lambda u}^\T\mu+\Qc_{t\to t',uu}u_t\big),
\end{align}
and the costate \eqref{eq:mr1} reads
\begin{align}
\pi_{t+1}={}&\Vc_{t+1\to t',xx}x_{t+1}+\Vc_{t+1\to t',\lambda x}^\T\mu+\vc_{t+1\to t',x}\notag\\
&+\Wc_{t+1\to t',x}^\T\nu_{t+1\to t'}.
\end{align}

\emph{Absorption.} With $T_\tau$, $\tilde W_{\tau\to t',x}$ and $\hat W_{\tau\to t',x}$ as in Section~\ref{sec:comb}, let
\begin{align}
\tilde w_{\tau\to t',0}&=\tilde T_\tau^\T\wc_{\tau\to t',0},\\
\hat w_{\tau\to t',0}&=\hat T_\tau^\T\wc_{\tau\to t',0},\\
S_\tau&=\tilde W_{\tau\to t',x}\Vc_{t\to\tau,\lambda\lambda}\tilde W_{\tau\to t',x}^\T.
\end{align}
The coefficients of \eqref{eq:numap} are
\begin{align}
N_{\tau,t\to\tau,x}&=-S_\tau^{-1}\tilde W_{\tau\to t',x}\Vc_{t\to\tau,\lambda x},\\
N_{\tau,t\to\tau,\lambda}&=-S_\tau^{-1}\tilde W_{\tau\to t',x}\Vc_{t\to\tau,\lambda\lambda},\\
n_{\tau,t\to\tau,0}&=-S_\tau^{-1}\big(\tilde W_{\tau\to t',x}\vc_{t\to\tau,\lambda}+\tilde w_{\tau\to t',0}\big),
\end{align}
and the absorbed tuple is
\begin{align}
\Vc'_{t\to\tau,xx}&=\Vc_{t\to\tau,xx}+\Vc_{t\to\tau,\lambda x}^\T\tilde W_{\tau\to t',x}^\T N_{\tau,t\to\tau,x},\\
\vc'_{t\to\tau,x}&=\vc_{t\to\tau,x}+\Vc_{t\to\tau,\lambda x}^\T\tilde W_{\tau\to t',x}^\T n_{\tau,t\to\tau,0},\\
\Vc'_{t\to\tau,\lambda x}&=\Vc_{t\to\tau,\lambda x}+\Vc_{t\to\tau,\lambda\lambda}\tilde W_{\tau\to t',x}^\T N_{\tau,t\to\tau,x},\\
\vc'_{t\to\tau,\lambda}&=\vc_{t\to\tau,\lambda}+\Vc_{t\to\tau,\lambda\lambda}\tilde W_{\tau\to t',x}^\T n_{\tau,t\to\tau,0},\\
\Vc'_{t\to\tau,\lambda\lambda}&=\Vc_{t\to\tau,\lambda\lambda}+\Vc_{t\to\tau,\lambda\lambda}\tilde W_{\tau\to t',x}^\T N_{\tau,t\to\tau,\lambda},\\
\Wc'_{t\to\tau,x}&=\begin{pmatrix}\Wc_{t\to\tau,x}\\\hat W_{\tau\to t',x}\Vc_{t\to\tau,\lambda x}\end{pmatrix},\\
\wc'_{t\to\tau,0}&=\begin{pmatrix}\wc_{t\to\tau,0}\\\hat W_{\tau\to t',x}\vc_{t\to\tau,\lambda}+\hat w_{\tau\to t',0}\end{pmatrix}.
\end{align}

\emph{Junction.} With
\begin{equation}
\Xi_\tau=\big(I-\Vc_{\tau\to t',xx}\Vc'_{t\to\tau,\lambda\lambda}\big)^{-1},
\end{equation}
the coefficients of \eqref{eq:lambdamap} are
\begin{align}
L_{\tau,t\to t',x}&=\Xi_\tau\Vc_{\tau\to t',xx}\Vc'_{t\to\tau,\lambda x},\\
L_{\tau,t\to t',\lambda}&=\Xi_\tau\Vc_{\tau\to t',\lambda x}^\T,\\
l_{\tau,t\to t',0}&=\Xi_\tau\big(\Vc_{\tau\to t',xx}\vc'_{t\to\tau,\lambda}+\vc_{\tau\to t',x}\big),\\
X_{\tau,t\to t',x}&=\Vc'_{t\to\tau,\lambda x}+\Vc'_{t\to\tau,\lambda\lambda}L_{\tau,t\to t',x},\\
X_{\tau,t\to t',\lambda}&=\Vc'_{t\to\tau,\lambda\lambda}L_{\tau,t\to t',\lambda},\\
x_{\tau,t\to t',0}&=\vc'_{t\to\tau,\lambda}+\Vc'_{t\to\tau,\lambda\lambda}l_{\tau,t\to t',0},
\end{align}
and the combined tuple is
\begin{align}
\Vc_{t\to t',xx}&=\Vc'_{t\to\tau,xx}+\Vc'^{\,\T}_{t\to\tau,\lambda x}L_{\tau,t\to t',x},\\
\vc_{t\to t',x}&=\vc'_{t\to\tau,x}+\Vc'^{\,\T}_{t\to\tau,\lambda x}l_{\tau,t\to t',0},\\
\Vc_{t\to t',\lambda x}&=\Vc_{\tau\to t',\lambda x}X_{\tau,t\to t',x},\\
\vc_{t\to t',\lambda}&=\vc_{\tau\to t',\lambda}+\Vc_{\tau\to t',\lambda x}x_{\tau,t\to t',0},\\
\Vc_{t\to t',\lambda\lambda}&=\Vc_{\tau\to t',\lambda\lambda}+\Vc_{\tau\to t',\lambda x}X_{\tau,t\to t',\lambda},\\
\Wc_{t\to t',x}&=\Wc'_{t\to\tau,x},\\
\wc_{t\to t',0}&=\wc'_{t\to\tau,0}.
\end{align}
For a conventional right tuple, $\Vc_{\tau\to t',xx}$, $\vc_{\tau\to t',x}$, $\Wc_{\tau\to t',x}$, $\wc_{\tau\to t',0}$ are replaced by
$V_{\tau,xx}$, $v_{\tau,x}$, $W_{\tau,x}$, $w_{\tau,0}$, all blocks in $\lambda_{t'}$ are dropped, and the maps carry the index $\tau,t$.

\emph{Down-sweep.} The composed maps of \eqref{eq:nutdown} are
\begin{align}
N_{\tau,t\to t',x}&=N_{\tau,t\to\tau,x}+N_{\tau,t\to\tau,\lambda}L_{\tau,t\to t',x},\\
N_{\tau,t\to t',\lambda}&=N_{\tau,t\to\tau,\lambda}L_{\tau,t\to t',\lambda},\\
n_{\tau,t\to t',0}&=n_{\tau,t\to\tau,0}+N_{\tau,t\to\tau,\lambda}l_{\tau,t\to t',0},\\
M_{\tau,t\to t',x}&=L_{\tau,t\to t',x}+\tilde W_{\tau\to t',x}^\T N_{\tau,t\to t',x},\\
M_{\tau,t\to t',\lambda}&=L_{\tau,t\to t',\lambda}+\tilde W_{\tau\to t',x}^\T N_{\tau,t\to t',\lambda},\\
m_{\tau,t\to t',0}&=l_{\tau,t\to t',0}+\tilde W_{\tau\to t',x}^\T n_{\tau,t\to t',0}.
\end{align}
The feedback policy of Corollary~\ref{cor:policy} at stage $t$ of leaf $i$, with $t'=t_{i+1}$, is $u^*_t(x_t)=U_{t,x}x_t+u_{t,0}$ with
\begin{align}
U_{t,x}&=U_{t\to t',x}+U_{t\to t',\lambda}M_{t',t,x},\\
u_{t,0}&=u_{t\to t',0}+U_{t\to t',\lambda}m_{t',t,0}.
\end{align}



\end{document}